\documentclass[conference]{IEEEtran}
\IEEEoverridecommandlockouts
\usepackage{cite}
\usepackage{amsmath,amssymb,amsfonts}
\allowdisplaybreaks
\usepackage{mathtools,comment}
\usepackage{algorithmic}
\usepackage{graphicx}
\usepackage{textcomp}
\usepackage{xcolor}
\let\originalleft\left
\let\originalright\right
\renewcommand{\left}{\mathopen{}\mathclose\bgroup\originalleft}
\renewcommand{\right}{\aftergroup\egroup\originalright}
\newcommand{\bE}{\mathbb{E}}

\newcommand{\bN}{\mathbb{N}}

\newcommand{\ust}{^{\star}}
\newcommand{\ub}{^{\beta}}
\newcommand{\bR}{\mathbb{R}}

\newcommand{\bP}{\mathbb{P}}

\newcommand{\cT}{\mathcal{T}}

\newcommand{\cF}{\mathcal{F}}

\newcommand{\Te}{\Tilde{e}}
\newcommand{\Tb}{\Tilde{b}}
\newcommand{\Tu}{\Tilde{u}}
\newcommand{\TV}{\Tilde{V}\ub}
\newcommand{\TQ}{\Tilde{Q}\ub}
\newcommand{\Tp}{\Tilde{p}}
\newcommand{\TVn}{\Tilde{V}\ub_n}
\newcommand{\TQn}{\Tilde{Q}\ub_{n+1}}
\newcommand{\itg}{\int_{\bR_+}}

\newcommand{\pe}{\psi(\Te_+)}

\newcommand{\ptot}{\psi(\Te_+,a\Te)}
\newcommand{\lf}{\left}
\newcommand{\rt}{\right}
\newcommand{\bl}{\bigl}
\newcommand{\br}{\bigr}
\newcommand{\enuma}{\renewcommand{\theenumi}{\alph{enumi}}}
\newcommand{\nal}[1]{\begin{align*}#1\end{align*}}
\newcommand{\al}[1]{\begin{align}#1\end{align}}
\newtheorem{assumption}{\textbf{Assumption}}

\newtheorem{definition}{Definition}[section]
\newtheorem{theorem}{Theorem}[section]
\newtheorem{proposition}{Proposition}[section]
\newtheorem{corollary}{Corollary}[theorem]
\newtheorem{lemma}[theorem]{Lemma}
\newtheorem{remark}{Remark}
\newcommand{\rom}[1]
    {\MakeUppercase{\romannumeral #1}}

\def\BibTeX{{\rm B\kern-.05em{\sc i\kern-.025em b}\kern-.08em
    T\kern-.1667em\lower.7ex\hbox{E}\kern-.125emX}}

\title{\LARGE \bf
Optimal Scheduling Policies for Remote Estimation of Autoregressive Markov Processes over Time-Correlated Fading Channel
}

\author{Manali Dutta and Rahul Singh
\thanks{$^{1}$Manali Dutta is with the Department of Electrical Communication Engineering, Indian Institute of Science, Bangalore, India
        {\tt\small manalidutta@iisc.ac.in}}%
\thanks{$^{2}$Rahul Singh is with the Department of Electrical Communication Engineering, Indian Institute of Science, Bangalore, India
        {\tt\small rahulsingh@iisc.ac.in}}%
}
\begin{document}

\maketitle

\begin{abstract}
We consider the problem of transmission scheduling for the remote estimation of a discrete-time autoregressive Markov process that is driven by white Gaussian noise. A sensor observes this process, and then decides to either encode the current state of this process into a data packet and attempts to transmit it to the estimator over an unreliable wireless channel modeled as a Gilbert-Elliott channel, or does not send any update. Each transmission attempt consumes $\lambda$ units of transmission power, and the remote estimator is assumed to be linear. The channel state is revealed only via the feedback (ACK\slash NACK) of a transmission, and hence the channel state is not revealed if no transmission occurs. The goal of the scheduler is to minimize the expected value of an infinite-horizon cumulative discounted cost, in which the instantaneous cost is composed of the following two quantities: (i)~squared estimation error, (ii) transmission power. We show that this problem can equivalently be posed as a partially observable Markov decision process (POMDP), in which the scheduler maintains a belief about the current state of the channel, and makes decisions on the basis of the current value of the estimation error, and the belief state.~We then show that the optimal policy is of threshold-type, i.e. for each value of the estimation error $e$, there is a threshold $b\ust(e)$ such that when the error is equal to $e$, then it is optimal to transmit only when the current belief state is greater than $b\ust(e)$.
\end{abstract}

\begin{IEEEkeywords}
Remote estimation, Gilbert-Elliott channel, partially observable Markov decision process (POMDP), threshold-type policy.
\end{IEEEkeywords}

\section{Introduction}
\subsection{Literature Overview}
In distributed networked control systems (NCS), several network nodes are connected through a communication network, which enables them to exchange information and collaborate to achieve a common goal~\cite{Hespanha2007survey,Antsaklis2007special}. In such a control architecture, decision-making is decentralized since each node can communicate only with its neighbors and makes decisions based on its own local information. Such systems have gained a widespread interest in recent years due to their ability to enable remote control and monitoring of physical systems. They are used in various fields, including industrial automation, robotics, and transportation systems.~Remote state estimation is one of the fundamental problems in NCS. Such a system is comprised of a sensor that observes an underlying random process, encodes its observations into data packets and then transmits it over a communication channel to a remote estimator that has a different location.~We will be exclusively interested in the case where the wireless medium is used for carrying out these transmissions. Since wireless devices are typically battery-operated, and transmissions consume energy, it is not efficient for the sensor to continually transmit the observations. Moreover, continuous transmissions can potentially flood the communication channel, leading to congestion and delays. Hence, in order to strike a balance between the communication cost and estimation error, sensors typically employ scheduling policies that make dynamic decisions regarding whether or not to send a packet, based on the information available with them.\par
In this work, we consider a remote estimator that attempts to estimate the state of a Markovian source in real-time. The works~\cite{Witsenhausen1979structure,Kaspi2012structure} investigate the structure of real-time encoders for a Markov source. Real-time encoders encode (quantize or compress) the output of a discrete time Markov source into a sequence of discrete variables, which in our case is a data packet. The encoder/sensor then transmits this sequence to the remote estimator. These operations are done in real-time. Several works have dealt with various aspects of the remote estimation problem.~For example,~\cite{Imer2010optimal,Xu2004optimal,Lipsa2011remote,Nayyar2013optimal,Chakravorty2014average}  consider the case when the communication channel is ideal, so that packet transmissions are always successful.~\cite{Imer2010optimal} studies the problem of estimating a discrete-time Markov process remotely, under the constraint that the sensor can access the wireless channel only a finite number of times. The sensor is restricted to the class of ``threshold type policies,'' i.e. those policies which transmit only when the current output at the source belongs to a particular set determined by the remaining number of channel uses, and the remaining number of decision instances.~The goal is to find an optimal estimator.~On the other hand,~\cite{Xu2004optimal} fixes the estimator to be ``Kalman-like'' and then optimizes over the choice of the scheduling policy for the sensor. It shows that the optimal scheduling decisions are solely a function of the current value of the state estimation error.~\cite{Lipsa2011remote}
does not impose any conditions on the structure of the scheduling policy or the estimator, and uses majorization theory~\cite{Marshall1979inequalities} in order to show that a threshold-type communication policy at the sensor, and a Kalman-like estimator are jointly optimal.~\cite{Nayyar2013optimal} shows that the structure of optimal communication and estimation policies derived in~\cite{Lipsa2011remote} continue to hold when additionally there are energy constraints on the transmitter.~\cite{Chakravorty2014average} also derives jointly optimal scheduling policy and estimator for the average cost problem by viewing it as a limiting case of
the discounted cost problem in the limit the discount factor approaches unity.~Communication policies which transmit only when the current value of the estimation error is greater than a threshold are also called event-triggered communication policies. Such communication policies are also analyzed and proposed for the case of error-free communication channel in \cite{Rabi2012adaptive,Molin2012iterative,Meng2012optimal,Han2015stochastic}.

Transmissions using the wireless medium are unreliable.~When communication is carried out over the wireless medium, the sensor needs to take this into account.~Various factors, such as environmental conditions, interference cause the state of the wireless channel to become ``bad,'' and this leads to packet losses. The works \cite{Chakravorty2016remote,Lipsa2009optimal,Chakravorty2017stochastic} study remote estimation problem for the case when the transmission occurs over wireless channels in which the transmitted packets are susceptible to losses; the packet losses are assumed to be i.i.d. across times. It is shown that the optimal transmission policy for such problems has a threshold structure with respect to the estimation error. A more realistic way to model the wireless fading channel is to model it using a finite state Markov chain\cite{Zhang1999finite}.~In this work, we model the unreliable wireless channel as a Gilbert-Elliott channel \cite{Gilbert1960capacity}, which is a Markovian channel in which the channel state can assume two values. At each discrete time, the channel is either in a good state and packet transmissions are successful, or it is in a bad state so that any attempted transmission fails.~The works~\cite{Ren2017infinite,Chakravorty2017structure,Chakravorty2019remote} study remote estimation over Markovian channels.~\cite{Ren2017infinite} derives transmission power control and remote estimation policies that are jointly optimal. It assumes that the channel state is instantaneously known to the sensor and estimator.~The problem is formulated as a partially observable Markov decision problem (POMDP), with a belief over the common information available with the sensor and the estimator, and the sensor makes decision based on the value of the belief state and the channel state. For the case when the power levels are discrete, it shows that the optimal transmission strategy is threshold-based with respect to the belief state. A model similar to \cite{Ren2017infinite} is considered in \cite{Chakravorty2017structure} and \cite{Chakravorty2019remote}, but with the difference that the channel state is known perfectly to the sensor with a delay of one unit. The optimality of a transmission policy that is of threshold-type with respect to the estimation error, is shown. However, obtaining a perfect knowledge of the channel state is difficult due to the complexity involved in measuring the characteristics of the communication channel. Factors such as physical limitations, cost, overheads and time-varying channel conditions, all contribute to the challenge of accurately measuring the channel state.

We address the problem of optimally scheduling transmissions to a linear estimator when the state of the channel is not completely observed by the sensor.~Thus, in our model the transmitter does not employ a probing mechanism to continually sense the channel state.~More concretely, the channel state is known to the transmitter only via acknowledgments that are sent by the estimator when it receives a packet upon a successful transmission. If there is no transmission attempt, then the current state is not known. The underlying process at the sensor which is being estimated is an autoregressive (AR) Markov process~\cite{Kumar2015stochastic,Chakravorty2019remote}, and our objective is to minimize the infinite-horizon cumulative expected discounted cost composed of (i) the squared estimation error and (ii) the transmission power. We formulate the problem of optimally scheduling transmissions as a POMDP. Motivated by the ease in implementation of policies that have a ``simple'' structure, we focus on studying structural properties of optimal transmission policies.~Since POMDPs are PSPACE hard~\cite{Papadimitriou1987complexity}, characterizing an optimal policy that has a simple structure reduces the search space. There is an extensive literature on structural results for POMDP~\cite{Lovejoy1987some,Albright1979structural,Krishnamurthy2007structured,Rieder1991structural,Grosfeld2007control,Ben2013partially,Krishnamurthy2009partially}. However, most of these works study POMDPs in which the state-space is just a simplex. However, in our work, the state-space is a Cartesian product of $\bR$ (the estimation error space) and $\bR_+$ (space in which the belief state lives). Moreover, in these works, the sensor makes decisions based on a noisy version of the underlying Markovian channel  state. Hence, we cannot apply the techniques used in these works. The work in  \cite{yao2022age} considers the problem of minimizing the long-run average age of information of a status update system under energy constraint and the channel is Gilbert-Elliott channel. They also consider the case where the channel state information is not available at the sensor while making decision, and is revealed  only via the feedback
(ACK/NACK) of a transmission. They formulate the problem as a POMDP with state-space consisting of completely observable states (age and time slot) which are discrete and can take only nonnegative values, and the belief state. In contrast, our work considers a continuous observable state-space (estimation error) that can take negative values. This complicates the analysis. 
\subsection{Contributions}
We consider the problem of designing optimal scheduling policies for a sensor that observes an autoregressive Markov process, and dynamically decides when to transmit these observations to a remote estimator, over an unreliable Markovian (Gilbert-Elliott)~\cite{Chakravorty2019remote,yao2022age,Chakravorty2017structure} wireless channel. The sensor gets to observe the channel state with a unit delay only when it attempts transmission and hence channel state is only partially observable. Our main contributions are as follows:
\begin{enumerate}
    \item We pose the problem faced by the sensor as a dynamic optimization problem that involves minimizing an infinite-horizon
    cumulative expected value of a discounted cost that consists of i) the squared estimation error and, ii) the transmission power.
    We show that this can be formulated as a Partially Observable MDP (POMDP)~\cite{Krishnamurthy2016partially} in which the state comprises of (a) the ``belief state,'' i.e. the conditional probability (conditioned on the information available with the sensor) that the channel state is good, (b) the current value of the estimation error. 
    \item Since our POMDP involves a one-stage cost function that is unbounded, it is not obvious that the value iteration algorithm~\cite{Puterman2014markov,Hernandez2012discrete} can be used to solve the POMDP. %
    We show that, under mild assumptions on the AR process and the Markovian transition probabilities of the channel~\eqref{assum:chan_plant}, the value iteration algorithm converges and yields an optimal policy.
    \item Since the analysis of the original POMDP is cumbersome, we introduce a certain ``folded POMDP," in which the state-space corresponding to the estimation error is $\bR_+$. We show the equivalence of the original POMDP and the folded POMDP, so that one can obtain an optimal policy and the value function for the original POMDP by solving the folded POMDP.~The concept of ``folding a Markov decision process (MDP)'' was introduced in \cite{Chakravorty2018sufficient}. However, since in our setup the channel state is not known by the sensor we cannot use the results of~\cite{Chakravorty2018sufficient}. 
    \item We then derive novel structural results for the POMDP by analyzing this folded POMDP. Specifically, we show that the optimal transmission strategy exhibits a threshold structure with respect to the belief state, and that there exists a threshold belief state such that it is optimal to transmit only when the current belief state is greater than that threshold.
\end{enumerate}

\section{Problem Formulation}
Consider a networked system comprising of a sensor and a remote estimator. The sensor observes an Auto Regression (AR)  process $x(t),t=0, 
 1,2,\ldots$ that evolves as follows,
$$
x(t+1) = a x(t) + w(t),
$$
where $a,x(t)\in\bR$, $w(t)$ is an i.i.d. Gaussian noise process that satisfies $w(t) \sim \mathcal{N}(0,1)$, with probability density function (pdf) given by $\varphi$, $\varphi(z) := \exp(-z^2\slash 2)\slash \sqrt{2\pi}$. Sensor encodes its observations into data packets, and transmits these to a remote estimator via an unreliable wireless channel.~We denote the state of the channel at time $t$ by $c(t)\in \{0,1\}$. $c(t)=0$ denotes that the channel is in ``bad state'' and any transmissions are unsuccessful, while $c(t)=1$ denotes that any packet transmitted at $t$ will be delivered to the estimator. The channel has \emph{memory}, and hence we assume that $\{c(t)\}_{t\in\bN}$ is a Markov process with parameters, 
\al{
    & p_{01} := \bP(c(t+1) =1 | c(t)=0), \label{eq:p01}\\
    & p_{11} := \bP(c(t+1) =1 | c(t)=1), \label{eq:p11}
}
where $p_{01},p_{11}>0$. Let $u(t)\in \{0,1\}$ denote the decision made by the sensor regarding whether ($u(t)=1$) or not ($u(t)=0$) a packet transmission should be attempted at time $t$.~We assume that each transmission attempt consumes $\lambda$ units of power\slash resource.~Let $y(t)$ denote the output of the channel, or the observation made by the estimator at time $t$, i.e.,
\al{
y(t) = 
\begin{cases}
    x(t) &\mbox{ if } c(t)u(t)=1,\\
    \Xi &\mbox{ if } c(t)u(t) = 0.
\end{cases}
}
where $y(t) = \Xi$ denotes that no packet was received, either because no transmission was carried out, or because the channel state was bad. Let $\hat{x}(t)$ denote the state of the estimator, or equivalently the point estimate made by the estimator. It evolves as follows,
\al{
\label{eq:estevolv}
\hat{x}(t+1) = 
\begin{cases}
    a \hat{x}(t) &\mbox{ if } y(t) = \Xi,\\
    y(t) &\mbox{ otherwise}. 
\end{cases}
}
Scheduler does not observe the channel state $c(t)$. However, if there is a successful transmission at $t$, then the estimator sends an acknowledgment to the sensor. Hence, if $u(t)=1$, then the channel state $c(t)$ at $t$ is known to the sensor at time $t+1$, or we say that upon transmitting a packet the scheduler gets to ``probe'' the channel. Letting 
$$
z(t) = 
\begin{cases}
    c(t) &\mbox{ if } u(t) = 1,\\
    \tilde{\Xi} &\mbox{ if } u(t) = 0,
\end{cases}
$$
where $\tilde{\Xi}$ denotes that no information about the current channel state was delivered, we have that the scheduler has access to $\{z(s)\}_{s=1}^{t-1}$ and also $\{u(s)\}_{s=1}^{t-1},\{x(s)\}_{s=1}^{t}$ while making the decision at time $t$. Here, $z(t) = \Xi$ denotes that the scheduler does not know the channel state when there is no transmission. Let 
\al{
    \label{def:error}
    e(t):= x(t) - \hat{x}(t),
}
denote the estimation error at $t$. The goal of the scheduler is to choose $\{u(t)\}_{t\in\bN}$ so as to minimize the expected value of cumulative estimation errors, as well as keep the cumulative transmission power utilized at minimal level.
~We let the instantaneous cost incurred by the system at time $t$ be $d(e(t)):=e(t)^2+\lambda u(t)$, and is a weighted sum of two quantities: 
(i) $e(t)^2$: square of the estimation error, (ii) communication cost $\lambda u(t)$. The goal of the scheduler at the sensor is to dynamically make the decisions $\{u(t)\}_{t\in\bN}$ so as to solve the following problem:
\al{
	\min_{\phi}~\bE_{\phi} \left(\sum_{t=0}^{\infty} \beta^{t}\left(e(t)^2+\lambda u(t) \right)\right),\label{def:obj}
}
where $\beta \in (0,1)$ is a discount factor, $\phi = \{\phi_t\}_{t\in\bN}$ is a measurable policy that for each time $t$ maps the history $\cF_t := \sigma( \{x(s),u(s),z(s)\}_{s=1}^{t-1}, x(t))$ to decision $u(t)$, and  $\bE_\phi$ denotes that the expectation is taken w.r.t the measure induced by the policy $\phi$.

\emph{Notation}: Let $\bN, \bR_{+}$ and $\bR_{-}$ denote the set of natural numbers, non-negative and non-positive real numbers, respectively.~$\delta_x(\cdot)$ is the delta function with unit mass at $x$.

\section{POMDP Formulation}
Note that while solving~\eqref{def:obj}, the channel state is not completely observed by the scheduler.~At each time $t$, it maintains an estimate of the current channel state, which is denoted by $b(t):= \bE(c(t)|\cF_t)$. This can be updated recursively using the ACK/NACK as follows,
\begin{equation} \label{def:beliefevolv}
	b(t+1)  = 
	\begin{cases}
		p_{11} & \text{if } u(t)=1,c(t)=1, \\
		p_{01} & \text{if } u(t)=1,c(t)=0, \\
        \cT(b(t)) & \text{if } u(t)=0,
	\end{cases}
\end{equation}
where for $x\in\bR$, we define $\cT(x) := x p_{11} + (1-x) p_{01}$. From~\eqref{eq:estevolv},~\eqref{def:error}, we have that the error process, $\{e(t)\}_{t \in \bN}$ evolves as follows,
\begin{equation} \label{def:error_evolve}
	e(t+1) =
	\begin{cases}
		ae(t)+w(t), & \text{if } u(t)c(t)=0 \\
		w(t), & \text{if } u(t)c(t)=1.
	\end{cases}
\end{equation}

For the purpose of solving~\eqref{def:obj}, we pose it as a POMDP~\cite{Krishnamurthy2016partially,Kumar2015stochastic}, in which the system state at $t$ is given by $(e(t),b(t))$,~where $e(t) \in \bR, b(t) \in [0,1]$ and $u(t) \in \{0,1\}$. The objective is to solve the following POMDP,
\al{
	\min_{\phi}~\bE_{\phi} \left(\sum_{t=0}^{\infty} \beta^{t}\left(e(t)^2+\lambda u(t) \right)\right),\label{def:pomdpobj}
}
where $b(t)$ and $e(t)$ evolve according to~\eqref{def:beliefevolv} and~\eqref{def:error_evolve}, respectively, the instantaneous cost incurred at time $t$ is given by,
\al{\label{def:instancost}
	d(e(t),b(t),u(t)) := e(t)^2 + \lambda u(t),
} 
and where at each time $t=0,1,2,\ldots,$ a scheduling policy $\phi$ takes action $u(t)$ on the basis of the operational history $\{(e(s),b(s))\}_{s=0}^{t}$.

We begin by discussing the POMDP transition model.~Let $p\left(e_+,b_+\mid e,b;u\right)$ denote the transition density function from the current state $(e,b)$ at time $t$ to the next state $(e_+,b_+)$ at time $t+1$ when action $u$ is taken at time $t$ for the above model. Consider the following two possibilities for $u$:

Case i) $u=0$: Then the state at the next step $(e_+,b_+)$ has the following density,
\al{
p(e_+,b_+ \mid e,b;0) = \exp(- (e_{+} - ae)^2\slash 2)\delta_{p_{11}b + p_{01}(1-b)}.\label{def:df_1}
}

Case ii) $u=1$: The density function of the resulting joint distribution of $(e_+,b_+)$ is as follows,
\al{
& p(e_+,b_+ \mid e,b;1) =b \exp(-e^2_+\slash 2)~\delta_{p_{11}}(b_+) \notag \\
& + (1-b) \exp(-(e_+ - ae  )^2\slash 2 )\delta_{p_{01}}(b_+).\label{def:df_2}
}

\subsection{Value Iteration}
We now show that under mild assumptions on the system parameters, value iteration algorithm can be used to solve the POMDP~\eqref{def:pomdpobj}. 
Value iteration algorithm is popularly used in order to solve MDPs. However, in order that we can use it to solve POMDP~\eqref{def:pomdpobj}, we need to verify whether our POMDP satisfies certain conditions~\cite[p.~46]{Hernandez2012discrete}. This is done below.
Define the $\beta$-discounted value function for the POMDP~\eqref{def:pomdpobj} as follows,
    \al{V\ub(e,b;\phi) := \bE_{\phi} \left(\sum_{t=0}^{\infty} \beta^{t}\left(e(t)^2+\lambda u(t) \right)\right), \label{Vphi}
    }
where $e \in \bR$ and $b \in [0,1]$. 
\begin{assumption}\label{assum:stabl}
    The Markovian channel probabilities and the system parameter $a$ satisfy the following condition
    \al{
    a^2(1-p_{01})<1.\label{assum:chan_plant}
    }
\end{assumption}
\begin{lemma} \label{lemma:assump}
Consider the POMDP~\eqref{def:pomdpobj}, and let Assumption~\ref{assum:stabl} hold. The following properties hold:
    \renewcommand{\labelenumi}{P\arabic{enumi}}
    \begin{enumerate}
    \item The one-stage cost function $e^2 + \lambda u $ is continuous, non-negative, and inf-compact on $(\bR \times [0,1] \times \{0,1\})$.
    \item The transition kernel from state $(e,b)$ at time $t$ to next state $(e_+,b_+)$ at time $t+1$  is strongly continuous for each action at time $t, u(t) = u \in 
    \{0,1\}$.
    \item There exists a policy $\phi$ such that $V\ub(e,b;\phi) < \infty$ for each $e \in \bR$ and $b \in [0,1]$.
    \end{enumerate}
\end{lemma}

The above result allows us to use value iteration. This is shown next. We begin by describing these iterations. Let $V\ub_n$ denote the value function at stage $n$ of the value iterations~\cite{Hernandez2012discrete}. We have for all $e \in \bR, b \in [0,1]$,
\begin{equation} 
    \label{eq:Vn}
 	V\ub_{n+1}(e,b)= \min_{u \in \{0,1\}} Q\ub_{n+1}(e,b;u),
 \end{equation}
where,
  \al{ 
            Q\ub_{n+1}(e,b;0) &:= e^2 + \beta \bE\left[V\ub_n(ae+w,\cT(b))\right]; \label{eq:Qn0}\\
            Q\ub_{n+1}(e,b;1) &:= e^2 + \lambda + \beta \bE\left[bV\ub_n(w,p_{11}) \right. \nonumber\\
            & \left. + (1-b)V\ub_n(ae+w,p_{01})\right] \label{eq:Qn1},
    }
    with,
    \begin{equation} \label{eq:V0}
            V\ub_0(e,b)=0.
        \end{equation}
Let $V\ub(e,b)$ denote the optimal total expected $\beta$-discounted cost function for the POMDP~\eqref{def:pomdpobj} 
, i.e.,
\al{
    \label{eq:optimalV}
    V\ub(e,b) : = \min_{\phi} V\ub(e,b;\phi)
}
The following proposition introduces the optimality equation for $V\ub$ and shows the convergence of value iteration method to $V\ub$.
\begin{proposition}\label{prop:vi}
Consider the POMDP~\eqref{def:pomdpobj} that satisfies Assumption~\ref{assum:stabl}. Then,
\begin{enumerate}
    \enuma
    \item Value iteration algorithm~\eqref{eq:Vn}-\eqref{eq:Qn1} converges to $V\ub$~\eqref{eq:optimalV}, i.e., 
    \begin{equation} \label{valfn}
        \lim_{n \rightarrow \infty} V\ub_n(e,b) = V\ub(e,b), \quad e \in \bR, b \in [0,1].
    \end{equation}
    \item Value function $V\ub$~\eqref{eq:optimalV} satisfies the following optimality equation,
        \begin{equation}    \label{eq:V}
            V\ub(e,b) = \min_{u \in \{0,1\}} Q\ub(e,b;u), \quad \forall e \in \bR, b \in [0,1],
        \end{equation}
        where,
        \begin{align} 
            Q\ub(e,b;0) &= e^2 + \beta \bE\left[V\ub(ae+w,\cT(b))\right], \label{eq:Q0}\\
            Q\ub(e,b;1) &= e^2 + \lambda + \beta \bE\left[bV\ub(w,p_{11}) \right. \nonumber\\
            & \left. + (1-b)V\ub(ae+w,p_{01})\right] \label{eq:Q1}.
        \end{align}
    \item There exists an optimal stationary deterministic policy that implements the minimizer of the right-hand side of~\eqref{eq:V} for each state $(e,b), e \in \bR, b \in [0,1]$ .    
\end{enumerate}   
\end{proposition}
\begin{IEEEproof}
    a) follows from~\cite[Lemma 4.2.8, pp. 49-59]{Hernandez2012discrete} since we have shown in Lemma~\ref{lemma:assump} that properties P1-P3 hold for POMDP~\eqref{def:pomdpobj}.~Similarly, b) and c) follow from~\cite[Theorem 4.2.3, pp. 46-47]{Hernandez2012discrete} since properties P1-P3 have been shown in Lemma~\ref{lemma:assump}.
    \end{IEEEproof}
\subsection{Folding the POMDP}
We will now derive some results for the POMDP~\eqref{def:pomdpobj} that allow us to ``fold it.'' This means that we construct an equivalent ``folded POMDP'' with state-space $\bR_+ \times [0,1]$, such that it suffices to study this POMDP in lieu of the original POMDP that has state-space $\bR \times [0,1]$. Specifically, the estimation error of the folded POMDP does not take negative values, in contrast to the original POMDP in which the estimation error takes both nonnegative and negative values.~Consequently, while analyzing the optimal policies, it is convenient to work with the folded POMDP rather than the original POMDP.~The work~\cite{Chakravorty2018sufficient} introduces the concept of a folded MDP. More specifically, for MDPs in which the state-space is $\bR$, it shows that under certain conditions on the transition probability kernel and instantaneous cost function, one can construct an equivalent MDP, called the ``folded MDP'' that has a state-space $\bR_{+}$ and is easier to study. Moreover, this also allows one to utilize an extensive theory on structural results for MDPs on $\bR_+$, or the set of natural numbers~\cite[Ch: 4,8]{Puterman2014markov}, in order to obtain structural results for an optimal policy for the original MDP.~However, the framework of~\cite{Chakravorty2018sufficient} cannot be used in order to study POMDPs.~Hence,~we now utilize the structure of POMDP~\eqref{def:pomdpobj} to introduce a ``folded POMDP.''
Before constructing the ``folded POMDP,'' we first show a structural property of the value function, $V\ub$ of the original POMDP~\eqref{def:pomdpobj}.
\begin{proposition}\label{prop1}
The functions $Q\ub(\cdot,b),V\ub(\cdot,b)$ for the POMDP~\eqref{def:pomdpobj} 
are even, i.e. we have $Q\ub(e,b;u)=Q\ub(|e|,b;u), V\ub(e,b)=V\ub(|e|,b)$ for all $b \in [0,1], u \in \{0,1\}$. 
\end{proposition}
\begin{IEEEproof}
We will use the properties of the iterates in~\eqref{eq:Vn}-\eqref{eq:V0} in order to prove this.~More specifically, since from~\eqref{valfn} we have that $\lim_{n \rightarrow \infty}V\ub_n(e,b) = V\ub(e,b)$, it suffices to show that $Q\ub_n(\cdot,b;u),V\ub_n(\cdot,b;u)$ are even for $n \in \bN$. This will then show that $Q\ub(\cdot,b;u),V\ub(\cdot,b)$ are also even. We will use induction in order to prove that $Q\ub_n,V\ub_n$ are even. Since $V\ub_0(e,b) = 0$ for all $e \in \bR$ and $b \in [0,1]$~\eqref{eq:V0}, $V\ub_0(\cdot,b)$ is even. Thus, the base case is true. Next, assume that the functions $V\ub_k(\cdot,b), b\in[0,1]$ are even for $k = 0,1,2,\ldots, n$. We will show that the functions $Q\ub_{n+1}(\cdot,b;u), b \in [0,1], u \in \{0,1\},$ are even. Consider the following two cases,

Case i): $u=0$. We have,
\al{
    & Q\ub_{n+1}(-e,b;0) = e^2 + \beta\bE\left[V\ub_n(-ae+w,\cT(b))\right]\label{eq:pr1}\\
    & = e^2 + \beta\int_{\bR}e^{-(e_+ + ae)^2 \slash 2}V\ub_n(e_+,\cT(b)) \,de_+ \label{eq:pr2}\\
    & = e^2 + \beta\int_{\bR} e^{-(-e' + ae)^2 \slash 2} V\ub_n(-e^{\prime},\cT(b)) \,de' \label{eq:pr3}\\
    & = e^2 + \beta\int_{\bR} e^{-(e' - ae)^2 \slash 2} V\ub_n(e^{\prime},\cT(b)) \,de' \label{eq:pr4}\\
    & = e^2 + \beta\bE\left[V\ub_n(ae+w,\cT(b))\right]\\
    &=Q\ub_{n+1}(e,b;u),
}
where \eqref{eq:pr1} follows from the definition of $Q\ub_n(e,b;0)$~\eqref{eq:Qn0}, while~\eqref{eq:pr2} follows from~\eqref{def:df_1}. The third equality~\eqref{eq:pr3} follows from a change of variables $e_+= -e^{\prime}$; while \eqref{eq:pr4} follows from our induction hypothesis that $V\ub_n(\cdot,b)$ is even. Thus, we have shown that $Q\ub_{n+1}(\cdot,b;0)$ is even.

Case ii): $u=1$. We have,
\al{
    & Q\ub_{n+1}(-e,b;1) = e^2 + \lambda \notag \\
    & + \beta \bE\left[bV\ub_n(w,p_{11}) + (1-b)V\ub_n(-ae+w,p_{01})\right] \label{eq:pr5}\\
    & = e^2 + \lambda + \beta b \int_{\bR} e^{-e_+^2 \slash 2}V\ub_n(e_+,p_{11}) \,de_+ \nonumber \\
    & + \beta (1-b) \int_{\bR} e^{-(e_+ + ae)^2 \slash 2}V_n\ub(e_+,p_{01}) \,de_+ \label{eq:pr6}\\
    & = e^2 + \lambda + \beta b \int_{\bR} e^{-e_+^2 \slash 2}V\ub_n(e_+,p_{11}) \,de_+ \nonumber \\
    & + \beta (1-b) \int_{\bR} e^{-(e_+ - ae)^2 \slash 2}V_n\ub(e_+,p_{01}) \,de_+ \label{eq:pr7}\\
    & = Q\ub_{n+1}(e,b;1), \label{eq:pr8}
}
where~\eqref{eq:pr5} follows from~\eqref{eq:Qn1}, while~\eqref{eq:pr6} follows from~\eqref{def:df_2}. The third equality follows from a change of variables and the induction hypothesis that $V\ub_n(\cdot,b)$ is even. This shows that $Q\ub_{n+1}(\cdot,b;1)$ is even.

So far we have shown that $Q\ub_{n+1}(\cdot,b;0),Q\ub_{n+1}(\cdot,b;1)$ are even.~Since $V\ub_{n+1}(\cdot,b)$ is the pointwise minimum of even functions $Q\ub_{n+1}(\cdot,b;0),Q\ub_{n+1}(\cdot,b;1)$~\eqref{eq:Vn}, it is even. The claim then follows by induction.
\end{IEEEproof}
We next dwell into the construction of ``folded POMDP.'' We use $\Tilde{\phi},\Tu,\Te$ and $\Tb$ to denote the policy, control, ``estimation error'' and ``belief state,'' respectively for the folded POMDP.\par
\emph{Folding the original POMDP}: We now construct a ``folded POMDP'' with the state-space $\mathbb{R}_+ \times[0,1]$; the error $\Te$ for this folded POMDP does not become negative. We will then show that on the set $\mathbb{R}_+ \times[0,1]$, the value functions and the optimal strategy of this folded POMDP is identical to that of the original POMDP. We then determine structural properties for the value function and optimal strategy for the folded POMDP, and then translate these properties back to the original POMDP~\eqref{def:pomdpobj}.\par
We begin with the definition of folded POMDP.
\begin{definition}[Folded POMDP]
    Given the original POMDP~\eqref{def:pomdpobj}, we define the ``folded POMDP'' on state-space $\bR_+ \times [0,1]$, control space $\{0,1\}$, and with the transition density function $\Tilde{p}$ defined as, 
    \al{
    \Tilde{p}\left(\Te_+,\Tb_+\mid \Te,\Tb;\Tu\right)= & ~p\left(\Te_+,\Tb_+\mid \Te,\Tb;\Tu\right) \notag \\
    & + p\left(-\Te_+,\Tb_+\mid \Te,\Tb;\Tu\right), \label{def:foldedPOMDP}
    }
    where $\Te,\Te_+ \in \mathbb{R}_+$, $\Tb,\Tb_+ \in [0,1]$ and $\Tu \in \{0,1\}$. The instantaneous cost, $d$ remains the same as in~\eqref{def:instancost}.
\end{definition}

We next show the equivalence of the original POMDP with state-space $\bR \times [0,1]$ and the folded POMDP with state-space $\bR_+ \times [0,1]$. We begin by discussing few properties of the folded POMDP.~Let $\Tilde{V}\ub(\Te,\Tb;\Tilde{\phi})$ and $\TV(\Te,\Tb)$ denote the $\beta$-discounted cost and $\beta$-discounted optimal value function, respectively, of the folded POMDP. These are analogous to~\eqref{def:pomdpobj} and~\eqref{eq:V}, respectively, of the original POMDP. We can show that the folded POMDP $(\mathbb{R}_+\times[0,1],\{0,1\},\Tilde{p},d)$ also satisfies the properties P1-P3 stated in Lemma \ref{lemma:assump}.~The proof is similar to that of Lemma \ref{lemma:assump}, which deals with the original POMDP~\eqref{def:pomdpobj}. Therefore, we can use value iteration to solve the folded POMDP also. Let $\TVn$ denote the iterates during stage $n$ of the value iteration algorithm~\cite{Hernandez2012discrete} when it is applied to solve the folded POMDP. We have the following for all $\Te \in \bR_+, \Tb \in [0,1]$,
\begin{equation} 
    \label{eq:Vntilde}
 	\TV_{n+1}\bl(\Te,\Tb\br)= \min_{\Tu \in \{0,1\}} \TQn\bl(\Te,\Tb;\Tu\br), 
 \end{equation}
    where, $\TQn\bl(\Te,\Tb;0\br)$ is as follows,
  \al{ 
            & \TQn\bl(\Te,\Tb;0\br) = \Te^2 + \beta \nonumber \\
            \times & \itg \Tp\left(\Te_+,\cT\bl(\Tb\br) \mid \Te,\Tb;0\right) \TVn\left(\Te_+,\cT\bl(\Tb\br)\right) \,d\Te_+, \label{eq1:Qn0tilde}\\
            = & ~\Te^2 + \beta \lf[\itg \!\!\!p\left(\Te_+,\cT\bl(\Tb\br) \mid \Te,\Tb;0\right) \TVn\left(\Te_+,\cT\bl(\Tb\br)\right) \,d\Te_+ \rt.\nonumber \\
            + & \lf. \!\itg \!\!\!p\left(-\Te_+,\cT\bl(\Tb\br) \mid \Te,\Tb;0\right) \TVn\left(\Te_+,\cT\bl(\Tb\br)\right) \,d\Te_+ \rt]  \label{eq2:Qn0tilde}\\
            = & ~\Te^2 + \beta \itg e^{-(\Te_+ - a\Te)^2 \slash 2} \TVn\left(\Te_+,\cT\bl(\Tb\br)\right) \,d\Te_+ \nonumber \\
            & + \beta \itg e^{-(\Te_+ + a\Te)^2 \slash 2} \TVn\left(\Te_+,\cT\bl(\Tb\br)\right) \,d\Te_+ , \label{eq:Qn0tilde}
        }
        where~\eqref{eq2:Qn0tilde} follows from the definition of folded POMDP~\eqref{def:foldedPOMDP}, while~\eqref{eq:Qn0tilde} follows from the definition of transition density in~\eqref{def:df_1}.

        While for $\Tu = 1$ we get,
        \al{
            &\TQn\bl(\Te,\Tb;1\br) \nonumber \\
            = & ~e^2 + \lambda + \nonumber \\
            + & \beta \lf[\Tb \itg \Tp\bl(\Te_+,p_{11} \mid \Te,\Tb;1\br) \TVn\left(\Te_+,p_{11}\right) \,d\Te_+ \rt. \nonumber\\
            + & \lf. (1-\Tb)\!\itg\!\!\! \Tp\bl(\Te_+,p_{01} \mid \Te,\Tb;1\br) \TVn\left(\Te_+,p_{01}\right) \,d\Te_+ \rt]  \label{eq1:Qn1tilde} \\
            = & ~\Te^2 + \lambda + \beta \Tb \itg 2e^{-(\Te_+)^2 \slash 2} \TVn\left(\Te_+,p_{11}\right) \,d\Te_+ \nonumber \\
            + & \beta (1-\Tb ) \itg e^{-(\Te_+ - a\Te)^2 \slash 2} \TVn\left(\Te_+,p_{01}\right) \,d\Te_+ \nonumber \\
            + & \beta (1-\Tb)\itg e^{-(\Te_+ + a\Te)^2 \slash 2} \TVn\left(\Te_+,p_{01}\right) \,d\Te_+ , \label{eq:Qn1tilde}
    }
    where~\eqref{eq:Qn1tilde} follows from~\eqref{def:foldedPOMDP} and~\eqref{def:df_2}. 
    
    The algorithm is initialized as follows,
    \begin{equation} \label{eq:V0tilde}
            \TV_0\bl(\Te,\Tb\br)=0, \quad \Te \in \bR_+, \Tb \in [0,1].
        \end{equation}
We have the following properties for the folded POMDP,~analogous to the results for the original POMDP shown in Proposition \ref{prop:vi}. These follow from~\cite[Theorem 4.2.3]{Hernandez2012discrete}.

\begin{enumerate}
    \enuma
    \item The value iteration algorithm with iterates $\TVn$ converges to $\TV$, i.e.
    \begin{equation} \label{valfntilde}
        \lim_{n \rightarrow \infty} \TVn\bl(\Te,\Tb\br) = \TV\bl(\Te,\Tb\br).
    \end{equation}
    \item The value function $\TV$ is the minimal bounded solution satisfying,
    \begin{equation}    \label{eq:Vtilde}
            \TV\bl(\Te,\Tb\br) = \min_{\Tu \in \{0,1\}} \TQ\bl(\Te,\Tb;\Tu\br), \quad \Te \in \bR_+,\Tb \in [0,1],
        \end{equation}
        where, 
        \al{ 
            & \TQ\bl(\Te,\Tb;0\br) \nonumber \\
            = & ~\Te^2 + \beta \notag \\
            \times & \itg \Tp\left(\Te_+,\cT\br(\Tb\br) \mid \Te,\Tb;0\right)\TV\left(\Te_+,\cT\bl(\Tb\br)\right) \,d\Te_+, \label{eq:Q0tilde}
            }
    and,
    \al{
            &\TQ\bl(\Te,\Tb;1\br) \nonumber \\
            = & ~\Te^2 + \lambda + \beta \left[\Tb \itg \!\!\!\Tp\bl(\Te_+,p_{11} \mid \Te,\Tb;1\br) \TVn\left(\Te_+,p_{11}\right) \,d\Te_+ \right. \nonumber\\
            + & \left. (1-\Tb)\!\itg \!\!\!\Tp\bl(\Te_+,p_{01}\! \mid\! \Te,\Tb;1\br) \TVn\left(\Te_+,p_{01}\right) \,d\Te_+ \right]\label{eq:Q1tilde}.
        }
    \item There exists an optimal stationary deterministic policy
    that implements the minimizer of the right-hand side of~\eqref{eq:Vtilde} in state $(\Te,\Tb)$, where $\Te \in \bR_+, \Tb \in [0, 1]$. 
\end{enumerate}
 
The following Proposition now shows the equivalence of folded POMDP with state-space $\bR_+ \times [0,1]$ and original POMDP with state-space $\bR \times [0,1]$.
\begin{proposition}
    \label{prop2}
The functions, $\Tilde{Q}\ub, \Tilde{V}\ub$ corresponding to the folded POMDP match with $Q\ub, V\ub$,~\eqref{eq:V}-\eqref{eq:Q1} of the original POMDP on $\mathbb{R}_+ \times[0,1]$, i.e., we have for all $e \in \mathbb{R}$, $b \in [0,1]$ and $u \in \{0,1\}$,
    \begin{align}\label{eq:prop2}
          Q\ub(e,b;u) = \Tilde{Q}\ub\left(|e|,b;u\right), V\ub(e,b)=\Tilde{V}\ub(|e|,b).
    \end{align}
\end{proposition}
\begin{IEEEproof}
    We will show these properties for the iterates obtained in value iteration, i.e. $Q\ub_n,V\ub_n,  n \in \bN$~\eqref{valfn}-\eqref{valfntilde} and $\Tilde{Q}\ub_n,\TVn$~~\eqref{eq:Vntilde}-\eqref{eq:V0tilde}. The result would then follow from~\eqref{eq:V} and~\eqref{eq:Vtilde}. We will use induction in order to prove this. We begin by analyzing the folded POMDP.
    
For $e \in \bR, b \in [0,1]$, we have $V\ub_0(e,b)=V\ub_0\left(|e|,b\right)=\Tilde{V}\ub_0\left(|e|,b\right) = 0$ by \eqref{eq:V0} and~\eqref{eq:V0tilde}, and hence, the base case holds. Next, assume that~\eqref{eq:prop2} holds for $k = 1,2,\ldots,n$. We will show that it also holds for time step $n+1$, and hence this will complete the induction. For $e\in\bR_+, b \in [0,1]$ and $u=0$, we have, 
    \al{
        & Q\ub_{n+1}(e,b;0) = e^2 + \beta\bE\left[V\ub_n(ae+w,\cT(b))\right] \label{eq1prop2}\\
        & = e^2 + \beta\int_{\bR}e^{-(e_+ - ae)^2 \slash 2}V\ub_n\left(e_+,\cT(b)\right) \,de_+ \label{eq2prop2}\\
        &= e^2 +  \int_{\bR_+}e^{-(e_+ - ae)^2 \slash 2}V\ub_n\left(e_+,\cT(b)\right) \,de_+ \nonumber \\
        & \qquad + \int_{\bR_-}e^{-(e_+ - ae)^2 \slash 2}V\ub_n\left(e_+,\cT(b)\right) \,de_+ \label{eq3prop2}\\
        &= e^2 + \int_{\bR_+}e^{-(e_+ - ae)^2 \slash 2}V\ub_n\left(e_+,\cT(b)\right) \,de_+ \nonumber \\
        & \qquad + \int_{\bR_+}e^{-(-e_+ - ae)^2 \slash 2}V\ub_n\left(-e_+,\cT(b)\right) \,de_+ \label{eq4prop2}\\
        &= e^2 + \int_{\bR_+}e^{-(e_+ - ae)^2 \slash 2}V\ub_n\left(e_+,\cT(b)\right) \,de_+ \nonumber \\
        & \qquad + \int_{\bR_+}e^{-(e_+ + ae)^2 \slash 2}V\ub_n\left(-e_+,\cT(b)\right) \,de_+ \label{eq5prop2}\\
        & = \Tilde{Q}\ub_{n+1}(e,b;0), \label{eq6prop2}
    }
where~\eqref{eq1prop2} follows from the definition of $Q\ub_{n+1}$~\eqref{eq:Qn0}. The second equality follows from~\eqref{def:df_1}.~\eqref{eq5prop2} holds since $V\ub_n$ is even~(Proposition \ref{prop1}), and from the induction hypothesis that for $e \in \bR_+, V\ub_n(e,b)=\Tilde{V}\ub_n(e,b)$. Finally,~\eqref{eq6prop2} follows from the definition of $\TQn$~\eqref{eq:Qn0tilde}. This shows $ \Tilde{Q}\ub_{n+1}(e,b;0) =  Q\ub_{n+1}(e,b;0)$. Similarly, we can show that $ \Tilde{Q}\ub_{n+1}(e,b;1) =  Q\ub_{n+1}(e,b;1)$.  From Proposition~\ref{prop1} we have $ Q\ub_{n+1}(e,b;u) =  Q\ub_{n+1}(-e,b;u)$, for $u=0$ and $u=1$, and hence we conclude $Q\ub_{n+1}(e,b;u)=\Tilde{Q}\ub_{n+1}(|e|,b;u)$ for all $e\in\bR$. Since the value function $V\ub_n,\tilde{V}\ub_n$~\eqref{eq:Vn},~\eqref{eq:Vntilde} are pointwise-minimum of the corresponding $Q$-functions, where the minimum is taken w.r.t. $u$, we also obtain $V\ub_{n+1}(e,b)=\tilde{V}\ub_{n+1}(|e|,b)$. This completes the induction step and also the proof.
\end{IEEEproof}

For ease of notation denote:
\nal{
& \psi(\Te_+) := e^{-{(\Te_+)^2} \slash 2 } \\
& \psi(\Te_+ - a\Te) := e^{-(\Te_+ - a\Te)^2 \slash 2}, \psi(\Te_+ + a\Te) := e^{-(\Te_+ + a\Te)^2 \slash 2}\\
& \psi(\Te_+,a\Te) := \psi(\Te_+ - a\Te) + \psi(\Te_+ + a\Te)
}
\section{Structural results of optimal policy}
Even though the folded POMDP is simpler and involves a reduced state-space, we are not able to utilize the existing works on structural results~\cite{Puterman2014markov,Krishnamurthy2016partially} in order to study structural properties of an optimal policy for POMDP~\eqref{def:pomdpobj}. Hence, we now utilize properties of the POMDP~\eqref{def:pomdpobj} to derive novel structural results for the folded POMDPs. Note that in departure with the existing works on structural results for POMDPs~\cite{Krishnamurthy2016partially}, the state-space of the folded POMDP is not just the simplex or $\bR_+$.\par
\begin{definition}[Threshold-type Policy]
We say that a scheduling policy for the folded POMDP $\Tilde{\phi}:\bR_{+}\times [0,1]\mapsto \{0,1\}$ is of threshold type if for each $\tilde{e}\in\bR_+$, there exists a threshold $b\ust(\tilde{e})$ such that when the current value of error is $\Te$, then it transmits only when the belief is greater than $b\ust(\Te)$.     
\end{definition}
The following is commonly assumed about the Gilbert-Elliott channels~\cite{yao2022age},~\cite{Abad2017channel},~\cite{Laourine2010betting},~and we will require this while analyzing properties of the optimal policy.
\begin{assumption}\label{assum:channel}
The Markovian channel parameters~\eqref{eq:p01},~\eqref{eq:p11} satisfy $p_{11} \geq p_{01}$. 
\end{assumption}

We now show that the optimal policy of the folded POMDP has a threshold-type structure.
\begin{theorem} \label{mainthm}
 Consider the folded POMDP $(\mathbb{R}_+\times[0,1],\{0,1\},\Tilde{p},d)$. Its value function $\Tilde{V}\ub$ satisfies the following properties:
	   \begin{enumerate}
    	\renewcommand{\theenumi}{(\Alph{enumi}}
	\item For each $\Tb$, the function $\Tilde{V}\ub\bl(\cdot,\Tb\br)$ is non-decreasing (with regards to $\Te$). 
	\item For each $\Te$, the function $\Tilde{V}\ub\left(\Te,\cdot\right)$ is non-increasing (with respect to $\Tb$). 
        \item For beliefs $x,y,z,\Tb$ such that $x \geq y$ and $z = \Tb x + (1-\Tb)y$, we have, \label{item:c}       
        \al{\label{eq:c}
            &(1-\Tb)\lambda + \Tb  \Tilde{V}\ub(\Te,x) \notag  \\
            &  + (1-\Tb) \Tilde{V}\ub(\Te,y)  
             \geq \Tilde{V}\ub(\Te,z) .
        }
        \item For each $\Te\in\bR_+$, there exists a threshold $\Tb\ust(\Te)$ such that it is optimal to transmit only when $\Tb \geq \Tb\ust(\Te)$. Thus, the optimal strategy corresponding to $\Tilde{V}\ub$ exhibits a threshold structure. 
    \end{enumerate}
\end{theorem}
\begin{IEEEproof}
    We will prove (A)-(D) for the iterates $\Tilde{V}\ub_n(\Te,\Tb), n\in \bN$ in~\eqref{eq:Vntilde}. We will show this via induction. The result would then follow from \eqref{valfntilde}, since we have $\lim_{n\to\infty} \Tilde{V}\ub_n(\Te,\Tb) = \Tilde{V}\ub(\Te,\Tb)$.
    
    Since $\Tilde{V}\ub_0\bl(\Te,\Tb\br)\equiv 0$~\eqref{eq:V0tilde}, (A)-(D) hold for $n=0$. Next, assume that (A)-(C) hold for $k=1,2,\ldots,n$.~The proof is divided into four steps.~We will firstly show that the threshold property~(D) holds for $k=n+1$, and then show (A)-(C) also hold for $k=n+1$. 

    \textbf{Step \rom{1}:} (D) holds for step $n+1$: We have $\Tilde{V}\ub_{n+1}\bl(\Te,\Tb\br)= \min_{\Tu \in \{0,1\}} \Tilde{Q}\ub_{n+1}\bl(\Te,
    \Tb;\Tu\br)$. Firstly, note that $\Tilde{Q}\ub_{n+1}\bl(\Te,\Tb;1\br)$ is a linear function of $\Tb$ by the definition of $\Tilde{Q}\ub_{n+1}$ in~\eqref{eq:Qn1tilde}. We will now show that $\Tilde{Q}\ub_{n+1}\bl(\Te,\Tb;0\br)$ is concave in $\Tb$.  Note that $\Tilde{V}\ub_n\bl(\Te,\Tb\br)$ is concave with respect to $\Tb$~\cite{Smallwood1973optimal}, so that for $\alpha \in [0,1]$ and beliefs $\Tb_1, \Tb_2\in [0,1]$, we have,
    \al{
        & \itg \ptot \left[\alpha \TVn\left(\Te_+,\cT\bl(\Tb_1\br)\right) \right. \notag \\
        & \left. + (1-\alpha)\TVn\left(\Te_+,\cT\bl(\Tb_2\br)\right)\right] \,d\Te_+ \nonumber\\
        & \geq \itg \ptot \TVn\left(\Te_+,\cT\bl(\alpha\Tb_1 + (1-\alpha)\Tb_2\br)\right) \,d\Te_+ \nonumber\\
        & = \!\!\itg \!\!\!\ptot \TVn\left(\Te_+, \alpha \cT\bl(\Tb_1\br)\!+\!(1-\alpha)\cT\bl(\Tb_2\br)\right) \!\,d\Te_+, \label{eq:cncv}
    }
    where the last equality follows from simple algebraic manipulations.  Concavity of $\Tilde{Q}\ub_{n+1}(\Te,\cdot;0)$ then follows from \eqref{eq:Qn0tilde} and \eqref{eq:cncv}.

Since $\lambda \geq 0$, from~\eqref{eq:Qn0tilde} and~\eqref{eq:Qn1tilde} we have that $\Tilde{Q}\ub_{n+1}(\Te,0;1) \geq \Tilde{Q}\ub_{n+1}(\Te,0;0)$. Now, consider the following two possible cases depending on the relationship between $\Tilde{Q}\ub_{n+1}(\Te,1;1)$ and $\Tilde{Q}\ub_{n+1}(\Te,1;0)$:
    
    Case i) $\Tilde{Q}\ub_{n+1}(\Te,1;1) < \Tilde{Q}\ub_{n+1}(\Te,1;0)$: then by the concavity of $\Tilde{Q}\ub_n\bl(\Te,\Tb;0\br)$ and linearity of $\Tilde{Q}\ub_{n+1}\bl(\Te,\Tb;1\br)$ in $\Tb$, there exists a unique point where the curves of $\Tilde{Q}\ub_{n+1}(\Te,1;1)$ and $\Tilde{Q}\ub_{n+1}(\Te,1;0)$ intersect. This intersection point corresponds to the threshold $\Tb\ust(\Te)$, i.e. during the $n+1$-th step of the iteration, it is optimal to transmit for belief values greater than this value.

    Case ii) $\Tilde{Q}\ub_{n+1}(\Te,1;1) \geq \Tilde{Q}\ub_{n+1}(\Te,1;0)$: we will show that for error value equal to $\Te$, it is optimal to not transmit for any value of $\Tb$. Specifically, we will prove that the curve of $\Tilde{Q}\ub_n(\Te,\cdot;1)$ always lies above the curve of $\Tilde{Q}\ub_n(\Te,\cdot;0)$, i.e. $\Tilde{Q}\ub_{n+1}\bl(\Te,\Tb;1\br) \geq \Tilde{Q}\ub_{n+1}\bl(\Te,\Tb;0\br)$ for all $\Tb\in [0,1]$. Now,
    \al{
         \Tilde{Q}\ub_{n+1}(\Te,1;1) \geq \Tilde{Q}\ub_{n+1}(\Te,1;0).
         }
Upon substituting~\eqref{eq:Qn0tilde},~\eqref{eq:Qn1tilde} into the above, we obtain,     
\al{ 
    & \lambda + \beta  \itg 2\pe \TVn(\Te_+,p_{11}) \,d\Te_+ \nonumber \\
    & \geq \beta  \itg \ptot \TVn(\Te_+,p_{11}) \,d\Te_+. \label{eq:step1}
    }
    Thus, we have
    \al{
    & \Tilde{Q}\ub\bl(\Te,\Tb;1\br)-\Tilde{Q}\ub\bl(\Te,\Tb;0\br) \notag\\
    & = \lambda + \beta \Tb \itg 2\pe \Tilde{V}\ub_n\left(\Te_+,p_{11}\right) \,d\Te_+ \nonumber \\
    & + \beta (1-\Tb) \itg \ptot \Tilde{V}\ub_n(\Te_+,p_{01}) \,d\Te_+ \notag\\
    & - \beta \itg \ptot \Tilde{V}\ub_n\left(\Te_+,\cT\bl(\Tb\br)\right) \,d\Te_+ \label{eq1:step1ii}\\
    & = \lambda + \beta \Tb \itg 2\pe \Tilde{V}\ub_n\left(\Te_+,p_{11}\right) \,d\Te_+ \nonumber \\
    & + \beta (1-\Tb) \itg \ptot \Tilde{V}\ub_n(\Te_+,p_{01}) \,d\Te_+ \notag \\
    & - \beta \itg \ptot \Tilde{V}\ub_n\left(\Te_+,\cT\bl(\Tb\br)\right) \,d\Te_+ \notag \\
    & + \Tb \beta \itg \pe \TVn(\Te_+,p_{11}) \,d\Te_+ \notag\\
    & - \Tb \beta \itg \pe \TVn(\Te_+,p_{11}) \,d\Te_+ + \Tb\lambda -\Tb\lambda \label{eq2:step1ii}\\
    & = \Tb\left[\lambda + \beta \itg 2\pe \Tilde{V}\ub_n(\Te_+,p_{11}) \,d\Te_+ \right. \notag \\
    & \left. - \beta \itg \ptot \Tilde{V}\ub_n(\Te_+,p_{11}) \,d\Te_+ \right] \notag\\
    & + \beta \left(\Tb \itg \ptot \Tilde{V}\ub_n(\Te_+,p_{11}) \,d\Te_+ \right. \notag \\
    & \left. + (1-\Tb) \itg \ptot \Tilde{V}\ub_n(\Te_+,p_{01}) \,d\Te_+ \right.\notag\\
    & \left.- \itg \ptot \Tilde{V}\ub_n\lf(\Te_+,\cT\bl(\Tb\br)\rt) \,d\Te_+ + (1-\Tb)\lambda \right) \notag\\
    & + (1-\beta)(1-\Tb)\lambda \label{eq3:step1ii}\\
    & \geq 0, \label{eq4:step1ii}
    } 
    where~\eqref{eq1:step1ii} follows from the definition of $\Tilde{Q}\ub_n$~\eqref{eq:Qn0tilde} and~\eqref{eq:Qn1tilde}, while~\eqref{eq4:step1ii} follows from (\ref{eq:step1}) and the induction hypothesis regarding property (C).
    
    \textbf{Step \rom{2}:} (A) holds for step $n+1$: Consider estimation errors $\Te,\Te' \in \bR_+$ satisfying $\Te^{\prime} > \Te$. We will show that $\Tilde{V}\ub_{n+1}\bl(\Te',\Tb\br) \geq \Tilde{V}\ub_{n+1}\bl(\Te,\Tb\br)$. From~\eqref{eq:Vntilde}, it suffices to show that for each value of control $\Tu \in \{0,1\}$ chosen for the state $\Te'$, there exists a control $\Tu' \in \{0,1\}$ under which the following holds, $\Tilde{Q}\ub_{n+1}\bl(\Te',\Tb;\Tu\br) \geq \Tilde{Q}\ub_{n+1}\bl(\Te,\Tb;\Tu'\br)$. We consider these two cases below separately.

    Case i): $\Tu=0$. 
    We have,
    \al{
    & \Tilde{Q}\ub_{n+1}(\Te',\Tb;0) \nonumber \\
    & = ({\Te'})^2 + \beta 
    \itg \psi(\Te_+,a\Te') \Tilde{V}\ub_n\lf(\Te_+,\cT\bl(\Tb\br)\rt) \,d\Te_+ \label{eq:19_0 } \\
    & \geq \Te^2 + \beta\itg \psi(\Te_+,a\Te) \Tilde{V}\ub_n\lf(\Te_+,\cT\bl(\Tb\br)\rt) \,d\Te_+  \label{eq:19}\\
    & = \Tilde{Q}\ub_{n+1}\bl(\Te,\Tb;0\br), \label{eq:20}
    }
    where~\eqref{eq:19_0 } follows from the definition of $\TQn$ in~\eqref{eq:Qn0tilde}, while~\eqref{eq:19} follows from Lemma \ref{lemmaA.1} in Appendix.

    Case ii): $\Tu=1$. We have,
    \al{
    & \Tilde{Q}\ub_{n+1}\bl(\Te',\Tb;1\br) \nonumber\\
    & = (\Te')^2 + \lambda + \beta \Tb \itg 2\pe \Tilde{V}\ub_n\left(\Te_+,p_{11}\right) \,d\Te_+ \nonumber \\
    & + \beta (1-\Tb) \itg \psi(\Te_+,a\Te') \Tilde{V}\ub_n(\Te_+,p_{01}) \,d\Te_+ 
    \label{eq1:Step2u=1} \\
    & \geq \Te^2 + \lambda + \beta \Tb \itg 2\pe \Tilde{V}\ub_n\left(\Te_+,p_{11}\right) \,d\Te_+ \nonumber \\
    & + \beta (1-\Tb) \itg \psi(\Te_+,a\Te) \Tilde{V}\ub_n(\Te_+,p_{01}) \,d\Te_+ 
    \label{eq2:Step2u=1}\\
    & = \Tilde{Q}\ub_{n+1}\bl(\Te,\Tb;1\br), \label{eq3:Step2u=1}
    }
    where~\eqref{eq1:Step2u=1} follows from~\eqref{eq:Qn1tilde}, while~\eqref{eq2:Step2u=1}  follows from Lemma \ref{lemmaA.1} in Appendix.
    
    \textbf{Step \rom{3}:} (B) holds for step $n+1$: Consider belief values $\Tb,\Tb'\in [0,1]$ satisfying $\Tb' \leq \Tb$. We will show that $\Tilde{V}\ub_{n+1}\bl(\Te,\Tb'\br) \geq \Tilde{V}\ub_{n+1}\bl(\Te,\Tb\br)$. To prove this, we will prove that for each value of control $\Tu$, we have $\Tilde{Q}\ub_{n+1}\bl(\Te,\Tb';\Tu\br)\ge \Tilde{Q}\ub\bl(\Te,\Tb;\Tu\br)$. Since $p_{11} \geq p_{01}$, we have $\cT\bl(\Tb'\br) \leq \cT\bl(\Tb\br)$. Consider the following two cases.
    
    Case i) $\Tu=0$: We have,
    \al{
    & \Tilde{Q}\ub_{n+1}\bl(\Te,\Tb';0\br) \nonumber \\
    & = \Te^2 + \beta 
    \itg \ptot \Tilde{V}\ub_n\lf(\Te_+,\cT\bl(\Tb'\br)\rt) \,d\Te_+ \label{eq1:Step3u=0 } \\
    & \geq \Te^2 + \beta\itg \ptot \Tilde{V}\ub_n\lf(\Te_+,\cT\bl(\Tb\br)\rt) \,d\Te_+  \label{eq2:Step3u=0}\\
    & = \Tilde{Q}\ub_{n+1}\bl(\Te,\Tb;0\br), \label{eq3:step3u=0}
    }
    where the first equality follows from~\eqref{eq:Qn0tilde}, while the inequality follows
    since (B) holds for $n$ by induction hypothesis.
    
    Case ii) $\Tu=1$: We have,
    \al{
     & \Tilde{Q}\ub_{n+1}\bl(\Te,\Tb';1\br) \nonumber \\
     = & ~\Te^2 + \lambda + \beta \itg \ptot \Tilde{V}\ub_n(\Te_+,p_{01})  \,d\Te_+ \nonumber \\
     + &  \beta \Tb' \itg \! \left(2\pe \Tilde{V}\ub_n(\Te_+,p_{11}) - \ptot \Tilde{V}\ub_n(\Te_+,p_{01})\right) \,d\Te_+ \nonumber \\
     \geq & ~\Te^2 + \lambda + \beta \itg \ptot \Tilde{V}\ub_n(\Te_+,p_{01})  \,d\Te_+  \nonumber \\
      + & \beta \Tb \itg  \left(2\pe \Tilde{V}\ub_n(\Te_+,p_{11}) \notag \right.\\
      & \qquad \qquad - \left.\ptot \Tilde{V}\ub_n(\Te_+,p_{01})\right) \,d\Te_+ \label{eq:25}\\
     = & ~Q\ub_{n+1}\bl(\Te,\Tb;1\br), \label{eq:26}
    }
where the first equality follows from the definition of $\TQn$ by~\eqref{eq:Qn1tilde}. By Lemma \ref{lemmaA.2} in Appendix, we have, $\itg \left(2\pe \Tilde{V}\ub_n(\Te_+,p_{11}) - \ptot \Tilde{V}\ub_n(\Te_+,p_{01})\right) \,d\Te_+ \leq 0$. Since, $\Tb' \leq \Tb$, the inequality~\eqref{eq:25} follows.

    \textbf{Step \rom{4}:} (C) holds for $n+1$: Now, since $x \geq y$, it follows from the threshold structure of policy which is optimal at stage $n+1$, proved in Step \rom{1}, that if the optimal action for state $(\Te,x)$ is to transmit, then the optimal action for state $(\Te,y)$ is also to transmit.~Thus, we have the following three possibilities while deciding optimal controls in states $(\Te,x)$ and~$(\Te,y)$, and we will separately show that this holds for all the cases:

    Case i) No transmission for both~$(\Te,x)$ and $(\Te,y)$: We have,
    \al{
    & (1-\Tb)\lambda + \Tb\TQn(\Te,x;0) + (1-\Tb)\TQn(\Te,y;0) \nonumber\\
    & = (1-\Tb)\lambda + \Tb\left[\Te^2 + \beta 
    \itg \ptot \Tilde{V}\ub_n\lf(\Te_+,\cT(x)\rt) \,d\Te_+\right] \nonumber\\
    & + (1-\Tb)\left[\Te^2 + \beta 
    \itg \ptot \Tilde{V}\ub_n\lf(\Te_+,\cT(y)\rt) \,d\Te_+\right] \nonumber \\
    & \geq \Te^2 + \beta 
    \itg \ptot \Tilde{V}\ub_n\lf(\Te_+,\cT(z)\rt) \,d\Te_+ \label{eq:30}\\
    & = \TQn(\Te,z;0) \label{eq:31}\\
    & \geq \Tilde{V}\ub_{n+1}(\Te,z), \label{eq:32}
    }
    where \eqref{eq:30} follows from the induction hypothesis on property \eqref{item:c}, while \eqref{eq:32} follows from \eqref{eq:Vntilde}.

    Case ii) Transmission for both the states $(\Te,x)$ and $(\Te,y)$: We have,
    \al{
    & (1-\Tb)\lambda + \Tb\TQn(\Te,x;1) + (1-\Tb)Q\ub_{n+1}(\Te,y;1) \nonumber\\
     & = (1-\Tb)\lambda + \Tb\lf(\Te^2 + \lambda + \beta x \itg 2\pe \Tilde{V}\ub_n\left(\Te_+,p_{11}\right) \,d\Te_+ \rt. \nonumber \\
    & \lf. + \beta (1-x) \itg \ptot \Tilde{V}\ub_n(\Te_+,p_{01}) \,d\Te_+\rt) \nonumber \\
    & + (1-\Tb)\lf(\Te^2 + \lambda + \beta y \itg 2\pe \Tilde{V}\ub_n\left(\Te_+,p_{11}\right) \,d\Te_+ \rt. \nonumber \\
    & \lf. + \beta (1-y) \itg \ptot \Tilde{V}\ub_n(\Te_+,p_{01}) \,d\Te_+\rt) \label{eq:33}\\
    & = (1-\Tb)\lambda + \Te^2 + \lambda \nonumber \\
    & + \beta \left[(\Tb x + (1-\Tb)y) \itg 2\pe \Tilde{V}\ub_n(\Te_+,p_{11})\right] \nonumber \\
    & + \beta \left[(1-\Tb x - (1-\Tb)y) \itg \ptot\Tilde{V}\ub_n(\Te_+,p_{01})\right] \label{eq:34} \\
    & = (1-\Tb)\lambda + \TQn(\Te,z;1) \label{eq:35} \\
    & \geq \TQn(\Te,z;1) \notag\\
    & \geq \Tilde{V}\ub_{n+1}(\Te,z), \label{eq:37}
    }
    where \eqref{eq:33} follows from definition \eqref{eq:Qn1tilde}, \eqref{eq:34} follows from some simple algebraic manipulations,~\eqref{eq:35} follows by the definition of $\TQn$ with $z=\Tb x + (1-\Tb)y$, and \eqref{eq:37} holds by \eqref{eq:Vntilde}.

    Case iii) Transmission for state $(\Te,x)$ and no transmission for state $(\Te,y)$: We have
    \al{
    & (1-\Tb)\lambda + \Tb \TQn(\Te,x;1) + (1-\Tb)\TQn(\Te,y;0) \nonumber \\
    & = \lambda + \Te^2 + \Tb\left[ \beta x \itg 2\pe \Tilde{V}\ub_n\left(\Te_+,p_{11}\right) \,d\Te_+ \rt. \nonumber \\
    & \lf. + \beta (1-x) \itg \ptot \Tilde{V}\ub_n(\Te_+,p_{01}) \,d\Te_+\right] \nonumber \\
    & + (1-\Tb)\left[\beta 
    \itg \ptot \Tilde{V}\ub_n\lf(\Te_+,\cT(y)\rt) \,d\Te_+\right] \label{eq:38} \\
    & \geq \lambda + \Te^2 + \Tb\beta\left[  x \itg 2\pe \Tilde{V}\ub_n\left(\Te_+,p_{11}\right) \,d\Te_+ \rt. \nonumber \\
    & \lf. + (1-x) \itg \ptot \Tilde{V}\ub_n(\Te_+,p_{01}) \,d\Te_+\right] \nonumber \\
    & + (1-\Tb)\beta \left[  y \itg \ptot \Tilde{V}\ub_n\left(\Te_+,p_{11}\right) \,d\Te_+ \rt. \nonumber \\
    & \lf. + (1-y) \itg \ptot \Tilde{V}\ub_n(\Te_+,p_{01}) \,d\Te_+\right] \label{eq:39} \\
    & \geq \lambda + \Te^2 + \beta \left[z \itg 2\pe \Tilde{V}\ub_n(\Te_+,p_{11}) \,d\Te_+ \rt. \notag \\
    & \lf. + (1-z)\itg \ptot \Tilde{V}\ub_n(\Te_+,p_{01})\right] \,d\Te_+ \label{eq:40} \\
    & = \TQn(\Te,z;1) \\
    & \geq \Tilde{V}\ub_{n+1}(\Te,z), \label{eq:42},
    }
    where \eqref{eq:38} follows from \eqref{eq:Qn0tilde} and \eqref{eq:Qn1tilde}. The inequality \eqref{eq:39} holds because $\Tilde{V}\ub_n$ is concave in $b$ from Step \rom{1} and \eqref{eq:40} follows from Remark \ref{remark1} after Lemma \ref{lemmaA.2} in the Appendix. Finally,~\eqref{eq:42} follows from \eqref{eq:Vntilde}.
\end{IEEEproof}
We now show that the POMDP~\eqref{def:pomdpobj} also admits an optimal policy that has a threshold structure.
\begin{corollary}
    \label{remark:folded-original}
     The original POMDP~\eqref{def:pomdpobj} satisfies the following properties:
     \begin{enumerate}
         \enuma
         \item The value function $V\ub$~\eqref{eq:V}  satisfies the properties (A)-(D) of Theorem \ref{mainthm}.
         \item The optimal strategy corresponding to $V\ub$ exhibits a threshold structure.
     \end{enumerate} 
\end{corollary}
\begin{IEEEproof}
    a) follows from Proposition \ref{prop2}. This is because we have shown that the folded POMDP with state space $\bR_+ \times [0,1]$ is equivalent to the original POMDP with state space $\bR \times [0,1]$.

    b) Let the optimal strategy corresponding to $V\ub$ and $\TV$ be $\phi\ust$ and $\Tilde{\phi}\ust$, respectively. Then, by property c) of Proposition \ref{prop:vi} we have,
    $\phi\ust(e,b) \in \mathop{\arg \min}\limits_{u \in \{0,1\}} Q\ub(e,b;u), e \in \bR, b \in [0,1]$. As a consequence of Proposition \ref{prop1}, we have that $\phi\ust(e,b)=\phi\ust(|e|,b)$, which means that the optimal strategy is even in $e$. Also, by Proposition \ref{prop2}, we have $Q\ub(e,b;u) = \TQ(|e|,b;u)$, which implies that $\phi\ust(e,b) = \Tilde{\phi}\ust(|e|,b)$. Now, since $\Tilde{\phi}\ust$ exhibits threshold structure by property (D) of Theorem \ref{mainthm}, it then follows that the optimal strategy $\phi\ust$ of the original POMDP~\eqref{def:pomdpobj} corresponding to $V\ub$ has a threshold structure such that $b\ust(e)=b\ust(|e|)=\Tb\ust(|e|)$.
\end{IEEEproof}

\section{Conclusion}
We considered a remote estimation problem in which the sensor observes an AR Markov process, and has to dynamically decide when to transmit updates to the estimator over a Gilbert-Elliott channel, so as to minimize a cumulative expected discounted cost that consists of estimation error and transmission power consumed. The sensor does not completely observe the channel, i.e. it obtains a delayed knowledge of the channel state only upon a transmission attempt. This problem can thus be posed as a POMDP, in which the decisions are solely a function of the current belief state and the estimation error. Since analyzing this POMDP is hard, we fold the POMDP, so that the ``error'' in the resulting POMDP remains positive. Consequently, we show an appealing structural result, namely that the optimal policy transmits only when the belief state is greater than a certain (error-dependent) threshold. This work can be extended in multiple directions. Firstly, a simple linear estimator is used, we would like to design an estimator and scheduler that are \emph{jointly} optimal. Secondly, the belief space is countably infinite, and hence the value iteration algorithm cannot be used in order to obtain the threshold values. We would like to obtain an efficient algorithm that would yield a good approximation to the optimal policy; one possibility could be to truncate the folded POMDP. We would also like to study a constrained remote estimation problem, in which there are constraints on average power consumption at the sensor. Finally, since the knowledge of AR process and Markovian channel parameters is not easy to obtain, we would like to design efficient learning algorithms which ``learn'' an estimator and scheduler that are jointly optimal asymptotically as $T\to\infty$.

\bibliographystyle{IEEEtran}
\bibliography{refs}

\appendices
\section{Proof of Lemma 3.1}
\label{app:A}
\begin{IEEEproof}
    P1) follows from the definition of the cost function $d$ and since our action set is finite. 
    
    P2) Let $P$ denote the transition kernel and suppose $\mu$ is the Lebesgue measure on $\bR$. Then, for any Borel measurable subset, B of $\bR$, we have by \cite[Example C.6]{Hernandez2012discrete},
\al{
    & P\bl((e_+,b_+) \in (B \times [0,1]) \mid e,b;u\br) \notag \\
    & = \sum_{b_+ \in [0,1]} \int_{B} p(e_+,b_+ \mid e,b;0) \mu(de_+) + \notag \\
    & + \sum_{b_+ \in [0,1]}\int_{B} p(e_+,b_+ \mid e,b;1) \mu(de_+) \label{eq1:trankernel}\\ 
    & = \sum_{b_+ \in [0,1]} \int_{B} p(e_+,b_+ \mid e,b;0) \,de_+ + \notag \\
    & + \sum_{b_+ \in [0,1]}\int_{B} p(e_+,b_+ \mid e,b;1) \,de_+, \label{eq2:trankernel}
}
    where~\eqref{eq1:trankernel} follows because the Lebesgue measure $\mu$ on $\bR$ is $\sigma$-finite and~\eqref{eq2:trankernel} follows because  $\mu(de_+) = de_+$.
    
    Then, $P$ is strongly continuous from the definition of $p$  \eqref{def:df_1}, \eqref{def:df_2}.
    
    P3) Consider the policy that transmits at every time step, i.e. $u(t)\equiv 1$. Consider the system starting in initial state $(e,b)$. The error at time $t$ can be written as follows,
\al{
e(t) = \left(\Pi_{m=1}^{t-1} a(m)\right) e + \sum_{s=1}^{t-1} \left(\Pi_{m=s}^{t-1} a(m)\right) w(s),
}
where 
\nal{
a(s):=
\begin{cases}
    a \mbox{ if } c(s) = 0\\
    0 \mbox{ if } c(s) = 1.
\end{cases}
}
Now, since $\{w(s)\}$ are i.i.d. and also independent of $\{a(s)\}$, we have that
\nal{
\bE~e(t)^2 = e^2 \bE\left(\Pi_{m=1}^{t-1} a(m)\right)^2 +  \bE\lf[\sum_{s=1}^{t-1} \left(\Pi_{m=s}^{t-1} a(m)\right)^2\rt].
}
We will now focus on $\bE \left(\Pi_{m=s}^{t-1} a(m)\right)^2$. Instead, consider $\bE\left\{ \left(\Pi_{m=s}^{t-1} a(m)\right)^2\Big| c(s)\right\}$. We have that $\Pi_{m=s}^{t-1} a(m)$ is equal to $0$ if $c(m)=1$ for atleast one $m\in \{s,s+1,\ldots,t-1\}$, and is equal to $a^{t-s}$ otherwise. The former occurs w.p. atleast $1-(1-p_{01})^{t-s}$ if $c(s)=0$, while when $c(s)=1$, this probability is atleast $p_{10}(1-p_{01})^{t-s-1}$. Upon using the law of total expectation, we obtain the following bound,
\al{
\bE \left(\Pi_{m=s}^{t-1} a(m)\right)^2 &\le (1-p_{01})^{t-s-1}a^{2(t-s)} \\
&= \frac{\left(a^2(1-p_{01})\right)^{t-s}}{1-p_{01}}. 
}
Since $a^2(1-p_{01})<1$, upon summing them above over $s$, we obtain the following, 
\al{
    \label{eq:assump1}
    \bE e(t)^2 \le \frac{1 + e^2}{\left(1-p_{01}\right)\left(1-a^2(1-p_{01})\right)}.
}

Since the cost per transmission is $\lambda$ units, from~\eqref{eq:assump1}, the cumulative discounted cost of the policy $u(t)\equiv 1$ is bounded by 
$$
\frac{1}{1-\beta} \left(\frac{1 + e^2}{\left(1-p_{01}\right)\left(1-a^2(1-p_{01})\right)} + \lambda \right).
$$
This completes the proof. 
\end{IEEEproof}
For ease of reference, we restate the notation here:
\nal{
& \psi(\Te_+) := e^{-{(\Te_+)^2} \slash 2 } \\
& \psi(\Te_+ - a\Te) := e^{-(\Te_+ - a\Te)^2 \slash 2}, \psi(\Te_+ + a\Te) := e^{-(\Te_+ + a\Te)^2 \slash 2}\\
& \psi(\Te_+,a\Te) := \psi(\Te_+ - a\Te) + \psi(\Te_+ + a\Te)
}

\begin{lemma} \label{lemmaA.1}
For $\Te',\Te \in \bR_+$ such that $\Te' \geq \Te$, the value iterates $\TVn$ corresponding to step $n$ in value iteration satisfies the following,
\al{
    \label{eq:lemmaA.1}
    & \itg \psi(\Te_+,a\Te') \Tilde{V}\ub_n\lf(\Te_+,\cT\bl(\Tb\br)\rt) \,d\Te_+ \nonumber \\
    & \geq \itg \ptot \Tilde{V}\ub_n\lf(\Te_+,\cT\bl(\Tb\br)\rt) \,d\Te_+.
} 
\end{lemma}
\begin{IEEEproof}
    To show~\eqref{eq:lemmaA.1} we need to divide it into two cases:

Case i) $a \geq 0$: Then, the term on the L.H.S of~\eqref{eq:lemmaA.1} is,
\al{
    & \int_{\bR_+} e^{-(\Te_+ - a\Te')^2 \slash 2} \Tilde{V}\ub_n\lf(\Te_+,\cT\bl(\Tb\br)\rt) \,d\Te_+ \nonumber\\
    & + \int_{\bR_+} e^{-(\Te_+ + a\Te')^2 \slash 2} \Tilde{V}\ub_n\lf(\Te_+,\cT\br(\Tb\br)\rt) \,d\Te_+. \label{eq3:a>0}
    }

For $\hat{e} > 0$ consider,
\al{
S(\hat{e},\Te):=\int_{\hat{e}}^{\infty} \left(e^{-(\Te_+ - a\Te)^2 \slash 2} + e^{-(\Te_+ + a\Te)^2 \slash 2}\right) \,d\Te_+, \label{def:S}
}
We will first show that~\eqref{def:S} is non-decreasing in $\Te$, i.e, for $\Te' \geq \Te$, we will have $S(\hat{e},\Te') \geq S(\hat{e},\Te)$.
Let $S_{\Te}(\hat{e},\Te)$ denote ${\partial S}\slash{\partial \Te}$. Then we have,
\al{
S_{\Te}(\hat{e},\Te) = a\left[e^{-(\hat{e} - a\Te)^2 \slash 2} - e^{-(\hat{e} + a\Te)^2 \slash 2}\right] \geq 0, \label{eq:dS,a>0}
}
where~\eqref{eq:dS,a>0} follows because $a \geq 0$ implies that $\left(e^{-(\hat{e} - a\Te)^2 \slash 2} - e^{-(\hat{e} + a\Te)^2 \slash 2}\right) \geq 0$.

Thus, we have $S(\hat{e},\Te') \geq S(\hat{e},\Te) $ with equality only if $\hat{e} = 0$. Now, since $\Tilde{V}\ub_n\bigl(\Te,\Tb\bigr)$ is non-decreasing in $\Te$ by induction hypothesis, we have by \cite[Lemma 4.7.2, p.106]{Puterman2014markov},
\nal{
& \int_{\bR_+} \left(e^{-(\Te_+ - a\Te')^2 \slash 2} + e^{-(\Te_+ + a\Te')^2 \slash 2}\right) \Tilde{V}\ub_n\lf(\Te_+,\cT\bl(\Tb\br)\rt) \,d\Te_+\\
& \geq \int_{\bR_+} \left(e^{-(\Te_+ - a\Te)^2 \slash 2} + e^{-(\Te_+ + a\Te)^2 \slash 2}\right) \Tilde{V}\ub_n\lf(e_+,\cT\bl(\Tb\br)\rt) \,d\Te_+,
}
which is exactly the claim.

Case ii) $a < 0$: The proof will follow along similar lines as
Case i) if we show that $S(\hat{e},\Te)$ defined in~\eqref{def:S} is non-decreasing in $\Te$. Now, we have,
\al{
S_{\Te}(\hat{e},\Te) = a\left[e^{-(\hat{e} - a\Te)^2 \slash 2} - e^{-(\hat{e} + a\Te)^2 \slash 2}\right] \geq 0, \label{eq:dS,a<0}
}
where~\eqref{eq:dS,a<0} follows because $a < 0$ implies that $\left(e^{-(\hat{e} - a\Te)^2 \slash 2} - e^{-(\hat{e} + a\Te)^2 \slash 2}\right) \leq 0$.
\end{IEEEproof}

\begin{lemma} \label{lemmaA.2}
		The value iterates $\TVn$ corresponding to step $n$ in value iteration satisfies the following,
		\begin{equation} \label{eq:lemmaA.2}
			\itg \left(2\pe \Tilde{V}\ub_n(\Te_+,p_{11}) - \ptot \Tilde{V}\ub_n(\Te_+,p_{01})\right) \,d\Te_+ \leq 0.
		\end{equation}
	\end{lemma}
\begin{IEEEproof}
     Firstly we note that,
     \al{
		& \itg \left(2\pe \Tilde{V}\ub_n(\Te_+,p_{11}) - \ptot \Tilde{V}\ub_n(\Te_+,p_{01})\right) \,d\Te_+ \nonumber\\
        \leq & \itg \bl(2\pe - \ptot \br) \Tilde{V}\ub_n(\Te_+,p_{01}) \,d\Te_+, \label{eq1:lemmaA3}
    }
    where~\eqref{eq1:lemmaA3} follows since $p_{11} \geq p_{01}$, and from the induction hypothesis that $\Tilde{V}\ub_n$ is non-increasing with respect to $b$.~For $\Hat{e} > 0$, consider,
    \nal{
    & \int_{\Hat{e}}^{\infty} \bl(2\pe - \ptot \br) \,d\Te_+ \\
    = & \int_{\Hat{e}}^{\infty}\!\! 2e^{-{\Te_+}^2 \slash 2} \,d\Te_+ - \!\int_{\Hat{e}}^{\infty}\!\! \left(e^{{-(\Te_+ -a\Te)^2 \!\slash 2}} + e^{-(\Te_+ + a\Te)^2 \!\slash 2}\right) \,d\Te_+
    }
     To show~\eqref{eq:lemmaA.2}, consider the following two cases:

    Case i) $a \geq 0$: 
    For $\Hat{e} > 0$, we further consider the following two cases:
    
    Case 1.i) $\Hat{e} \geq a\Te$: We have,
    \al{
        & \int_{\Hat{e}}^{\infty}\!\! 2e^{-{\Te_+}^2 \slash 2} \,d\Te_+ - \!\int_{\Hat{e}}^{\infty}\!\! \left(e^{{-(\Te_+ -a\Te)^2 \!\slash 2}} + e^{-(\Te_+ + a\Te)^2 \!\slash 2}\right) \,d\Te_+ \notag \\
        = & \int_{\Hat{e}}^{\infty} 2e^{-{\Te_+}^2 \slash 2} \,d\Te_+ - \int_{\Hat{e}-a\Te}^{\infty} e^{-{\Te_+}^2 \slash 2} \,d\Te_+ - 
        \int_{\Hat{e}+a\Te}^{\infty} e^{-{\Te_+}^2 \slash 2} \,d\Te_+ \notag \\
        = & \int_{\Hat{e}}^{\Hat{e}+a\Te} e^{-{\Te_+}^2 \slash 2} \,d\Te_+ - \int_{\Hat{e}-a\Te}^{\Hat{e}} e^{-{\Te_+}^2 \slash 2} \,d\Te_+ \notag \\
        \leq & ~0. \label{eq1:lemma2}
    }

    Case 2.i) $\Hat{e} < a\Te$: We have,
     \al{
        & \int_{\Hat{e}}^{\infty}\!\! 2e^{-{\Te_+}^2 \slash 2} \,d\Te_+ - \int_{\Hat{e}}^{\infty} \!\left(e^{{-(\Te_+ -a\Te)^2 \slash 2}} + e^{-(\Te_+ + a\Te)^2 \slash 2}\right) \,d\Te_+ \notag \\
        & = \int_{\Hat{e}}^{\Hat{e}+a\Te} e^{-{\Te_+}^2 \slash 2} \,d\Te_+ \nonumber \\
        & - \lf (\int_{0}^{\Hat{e}} e^{-{\Te_+}^2 \slash 2} \,d\Te_+ + \int_{0}^{a\Te -\Hat{e}} e^{-{\Te_+}^2 \slash 2} \,d\Te_+\rt) \nonumber \notag\\
        & \leq 0. \label{eq2:lemma2}
        }
    
    It then follows from inequalities~\eqref{eq1:lemma2} and~\eqref{eq2:lemma2} that for $\Hat{e} > 0$, we have,
		\al{
			& \int_{\Hat{e}}^{\infty} 2e^{-{\Te_+}^2 \slash 2} \,d\Te_+\nonumber \\
			& \leq \int_{\Hat{e}}^{\infty} \left(e^{{-(\Te_+ -a\Te)^2 \slash 2}} + e^{-(\Te_+ + a\Te)^2 \slash 2}\right) \,d\Te_+, \label{eq4:lemmaA3}
		}
    with equality holding in \eqref{eq4:lemmaA3} only when $\Hat{e} = 0$.
    Since, $\Tilde{V}\bl(\Te,\Tb\br)$ in non-decreasing in $\Te$ for $\Te \in \bR_+$ shown in Step (\rom{2}) of Theorem~\ref{mainthm}, by \eqref{eq4:lemmaA3}, and \cite[Lemma 4.7.2, p.106]{Puterman2014markov}, we have,
    \nal{
        & \int_{\bR_+} 2e^{-{\Te_+}^2 \slash 2} \Tilde{V}_n(\Te_+,p_{01}) \,d\Te_+ \nonumber \\
		& \leq \int_{\bR_+} \left(e^{{-(\Te_+ -a\Te)^2 \slash 2}} + e^{-(\Te_+ + a\Te)^2 \slash 2}\right) \Tilde{V}_n(\Te_+,p_{01}) \,d\Te_+.
    }
    Thus, the claim holds follows  from~\eqref{eq1:lemmaA3}.
    
    Case ii) $a < 0$: The proof will follow along similar lines as Case i) if we show that~\eqref{eq4:lemmaA3} holds.
    
    For $\Hat{e} > 0$, we consider the following two cases,

    Case 1.ii) $|\Hat{e}| \geq |a\Te|$: We have,
    \al{
        & \int_{\Hat{e}}^{\infty}\!\! 2e^{-{\Te_+}^2 \slash 2} \,d\Te_+ - \!\int_{\Hat{e}}^{\infty}\!\! \left(e^{{-(\Te_+ -a\Te)^2 \!\slash 2}} + e^{-(\Te_+ + a\Te)^2 \!\slash 2}\right) \,d\Te_+ \notag \\
        & = \int_{\Tilde{e}}^{\Tilde{e}-a\Te} e^{-{e_+}^2 \slash 2} \,de_+ - \int_{\Tilde{e}+ae}^{\Tilde{e}} e^{-{e_+}^2 \slash 2} \,de_+ \notag \\
        & \leq 0. \label{eq1:caseii.1}
    }

    Case 2.ii) $|\Hat{e}| < |a\Te|$: We have,
    \al{
        & \int_{\Hat{e}}^{\infty}\!\! 2e^{-{\Te_+}^2 \slash 2} \,d\Te_+ - \int_{\Hat{e}}^{\infty} \!\left(e^{{-(\Te_+ -a\Te)^2 \slash 2}} + e^{-(\Te_+ + a\Te)^2 \slash 2}\right) \,d\Te_+ \notag \\
        & = \int_{\Hat{e}}^{\Hat{e}-a\Te} e^{-{\Te_+}^2 \slash 2} \,d\Te_+ \nonumber \\
        & - \lf (\int_{0}^{\Hat{e}} e^{-{\Te_+}^2 \slash 2} \,d\Te_+ + \int_{0}^{-a\Te -\Hat{e}} e^{-{\Te_+}^2 \slash 2} \,d\Te_+\rt) \nonumber \notag \\
        & \leq 0. \label{eq2:caseii.2}
        }
    Thus, from inequalities~\eqref{eq1:caseii.1} and~\eqref{eq2:caseii.2} we have that,
		\al{
			& \int_{\Hat{e}}^{\infty} 2e^{-{\Te_+}^2 \slash 2} \,d\Te_+\nonumber \\
			& \leq \int_{\Hat{e}}^{\infty} \left(e^{{-(\Te_+ -a\Te)^2 \slash 2}} + e^{-(\Te_+ + a\Te)^2 \slash 2}\right) \,d\Te_+, 
		 \label{eq:l3}
		}
    with equality holding in \eqref{eq:l3} only when $\Hat{e} = 0$.
    
\end{IEEEproof}
\begin{remark} \label{remark1}
    Note that from~\eqref{eq1:lemmaA3}, the result of Lemma \ref{lemmaA.2} holds for any $\Tb \in [0,1]$, i.e.,
    \nal{
    & \itg \left(2\pe \Tilde{V}\ub_n(\Te_+,\Tb) - \ptot \Tilde{V}\ub_n(\Te_+,\Tb)\right) \,d\Te_+ \\
    & \leq 0.
    }
\end{remark}

\end{document}